\documentclass[10pt]{amsart}
\usepackage{amssymb,amsthm} 
\usepackage{pxfonts}   
\usepackage[dvips]{graphicx}
\usepackage{psfrag}    

\newtheorem{thm}{Theorem}
\newtheorem{corol}{Corollary}
\newtheorem{prop}{Proposition}
\newtheorem{lemma}{Lemma}
\theoremstyle{definition}
\newtheorem{defn}{Definition}

\theoremstyle{remark}
\newtheorem*{ack}{Acknowledgement}

\def\R{\mathbb{R}}

\def\Z{\mathbb{Z}}
\def\N{\mathbb{N}}
\def\T{\mathbb{T}}

\def\cE{\mathcal{E}}

\def\cH{\mathcal{H}}
\def\cI{\mathcal{I}}

\def\cP{\mathcal{P}}
\def\cR{\mathcal{R}}
\def\cU{\mathcal{U}}

\def\order{\mathrm{order}}
\renewcommand{\epsilon}{\varepsilon}
\renewcommand{\setminus}{\smallsetminus}
\renewcommand{\emptyset}{\varnothing}
\def\Elip{{\cE_0^\mathrm{lip}}}
\def\Epl{{\cE_0^\mathrm{pl}}}
\DeclareMathOperator*{\essinf}{ess \, inf}
\DeclareMathOperator*{\esssup}{ess \, sup}
\def\Dist{{\mathrm{Dist}}}
\def\dist{{d_\mathrm{lip}}}

\begin{document}

\title[Expanding maps without aci$\sigma$]{Generic expanding maps
without absolutely continuous invariant $\sigma$-finite measure}


\begin{abstract}

We show that a $C^1$-generic expanding
map of the circle has no absolutely continuous invariant
$\sigma$-finite measure.

\end{abstract}

\author[A.~Avila]{Artur Avila}
\address{CNRS UMR 7599,
Laboratoire de Probabilit\'es et Mod\`eles al\'eatoires.
Universit\'e Pierre et Marie Curie--Bo\^\i te courrier 188.
75252--Paris Cedex 05, France}
\urladdr{www.proba.jussieu.fr/pageperso/artur/}
\email{artur@ccr.jussieu.fr}

\author[J.~Bochi]{Jairo Bochi}
\address{Instituto de Matem\'atica -- UFRGS -- Porto Alegre, Brazil}
\urladdr{www.mat.ufrgs.br/$\sim$jairo}
\email{jairo@mat.ufrgs.br}

\date{August 2, 2006}

\maketitle

\section{Introduction}

If $f$ is a measurable transformation of a Lebesgue measure space $(X,\mathcal{A}, \lambda)$ to itself,
that does not preserve the measure $\lambda$,
one can study the invariant measures of $f$ and compare them to $\lambda$.
A especially interesting case is when $f$ is non-singular with respect to $\lambda$
(in the sense that $\lambda(A)=0$ iff $\lambda(f^{-1}(A))=0$),
but nevertheless there exist no $\sigma$-finite invariant measure which is absolutely
continuous with respect to $\lambda$.
Such maps $f$ are called of \emph{type III} (with respect to the measure).
Their existence was conjectured by Halmos~\cite{Halmos}
and established by Ornstein~\cite{O}.
Other examples were given later; 
let us cite a few (when not specified, the relevant measure is Riemannian):
\begin{itemize}
\item piecewise linear homeomorphisms of the circle, by Herman~\cite{Herman};
\item $C^\infty$-diffeomorphisms of the circle, by Katznelson~\cite{K};
\item a $C^\infty$ non-invertible map of the $2$-torus, by Hawkins and Silva~\cite{HSi};
\item the full shift on $2$ symbols, with respect to some product measure,
by Hamachi~\cite{Hamachi}.
\item a $C^1$ expanding map of the circle (constructed using Hamachi's
example), by Bruin and Hawkins~\cite{BruinHawkins}.
\end{itemize}
Recall that $C^{1+\alpha}$ expanding maps have
absolutely continuous invariant probability measures, so the regularity of
the example of Bruin and Hawkins is essentially sharp.

The question of whether the absence of aci$\sigma$ is actually a generic
(in the usual topological sense) phenomenon
for $C^1$ expanding maps of the circle
seems to have been first posed by Quas~\cite {Quas}.
Later investigations~\cite {CampbellQuas} indicated that the known methods
failed to decide the question either way.  It was also known that
$C^1$-generic maps do behave ``pathologically'' in some
respects (they have no absolutely continuous invariant {\it probability}
measure~\cite {Quas}), but not in others
(they are ergodic and conservative with respect to Lebesgue measure~\cite
{Quas}, and they possess a unique physical measure~\cite {CampbellQuas}).

In this paper we show that the type III property is indeed
$C^1$-generic for expanding maps of the circle.



\smallskip

We also mention that the non-existence of \emph{finite} 
invariant measures that are absolutely continuous
with respect to Riemannian measure was shown to be a generic property among 
(expanding or not) $C^1$ maps of compact manifolds of any dimension -- see~\cite{nos}.

Of course, it is natural to ask whether 
the result of the present paper is still true for expanding maps on higher dimension.
It is not clear whether our methods can be extended.

Concerning non-necessarily expanding maps of a compact manifold,
there are $C^1$-open sets of transformations that
do have some absolutely continuous $\sigma$-finite invariant measure
(maps with a sink, for instance).
In this regard, we ask whether a $C^1$-generic map
has no absolutely continuous $\sigma$-finite invariant measure
which is conservative (all w.r.t.~Riemannian measure).  We will show
that this is true at least for one-dimensional maps, see corollary \ref{besta}.

\smallskip

Let us now give the precise statements.

Let $\T^1 = \R / \Z$ be the circle.
Let $\cE^1$ be the set of all  $C^1$ maps $f:\T^1 \to \T^1$ which are expanding, i.e.,
\begin{equation} \label{e.c}
\exists c>0, \ \exists \lambda>1 \text{ s.t. } |(f^n)'(x)| > c\lambda^n \ \forall x \in\T^1, \ \forall n \in \N \, .
\end{equation}

We endow the set $\cE^1$ with the $C^1$ topology.
Let $m$ denote the Lebesgue measure on $\T^1$ normalized so that $m(\T^1)=1$.
We say that a $\sigma$-finite measure on $\T^1$ is an \emph{aci$\sigma$}
for a map $f:\T^1 \to \T^1$
if it is absolutely continuous with respect to $m$ and it is $f$-invariant.
Our main result is:

\begin{thm}\label{t.main}
There exists a residual set $\cR \subset \cE^1$ such that
if $f \in \cR$ then $f$ has no aci$\sigma$.
\end{thm}

\begin{corol} \label {besta}

Let $X$ be either the circle $\T^1$ or the compact interval $[0,1]$.
There exists a residual set $\cR'$
of the space of $C^1(X,X)$ such
that if $f \in \cR'$ then $f$ has no conservative aci$\sigma$.

\end{corol}

\begin{proof}

Hyperbolic maps form an open and dense subset $\cH$ of $C^1(X,X)$ by
\cite {J}. (See also \cite {KSS} for the recent extension to higher regularity.)
The map that associates to $f \in \cH$
its non-wandering set $\Omega(f)$ is upper semi-continuous in the Hausdorff
topology.  Moreover, if $f \in \cH \cap C^2(X,X)$ then $m(\Omega(f))=0$
unless $X=\T^1$ and $f$ is expanding, see \cite {M}.  It follows that generically,
either $m(\Omega(f))=0$ or $f \in \cE^1$.  In the first case, $f$
cannot have a conservative aci$\sigma$ (since any conservative measure
must be supported on the non-wandering set).  In the second case, generically
there is no aci$\sigma$ at all, by theorem \ref{t.main}.
\end{proof}

\section{Some preliminaries}

\subsection{A reduction}

Let $\cE^1_0$ be the
subset of $\cE^1$ consisting of maps $f$ satisfying $f(0)=0$ and
such that~\eqref{e.c} holds with $c=1$.  We will actually prove:

\begin{thm}\label{t.reduced}
There exists a residual set $\cR_0 \subset \cE^1_0$ such that if $f \in \cR_0$ then
$f$ has no aci$\sigma$.
\end{thm}

Let us show that theorem~\ref{t.reduced} implies theorem~\ref{t.main}. 
Let $f \in \cE^1$, and let $p_f$ be a fixed point of $f$.  Then for $\tilde f$ in
a small open neighborhood $\cU$ of $f$ in $\cE^1$,
there exists a unique fixed point $p_{\tilde f}$
of $\tilde f$ near $p_f$, moreover this fixed point depends continuously on $\tilde f$.
Let $n \geq 1$ be such that $|(f^n)'(x)|>1$ for every $x \in \T^1$.
Define 
$$
h(x)=\sum_{k=0}^{n-1} |(f^k)'(x)| \quad \text{so that} \quad
\int_{\T^1} h(x) dx=\frac {|d_f|^n-1} {|d_f|-1},
$$ 
where $d_f$ is the degree of $f$.  Let $H_{\tilde f}$ be the orientation preserving
$C^1$ diffeomorphism of $\T^1$ such that $H_{\tilde f}(p_{\tilde f})=0$ and
$H_{\tilde f}'(x)=\frac {|d_f|-1} {|d_f|^n-1} h(x)$.
Then $H_{\tilde f}$
depends continuously on $\tilde f \in \cE^1$.
Define a map $g=H_f \circ f \circ H^{-1}_f$.
For any $x \in \T^1$, writing $y = H_f^{-1}(x)$, we have 
$$
|g'(x)| = \frac{h(f(y)) \cdot |f'(y)|}{h(y)}
= \frac{h(y) +|(f^n)'(y)| - 1}{h(y)} > 1 .
$$
Hence $g \in \cE^1_0$.
Shrinking $\cU$ to a smaller open neighborhood of $f$ if necessary, we see that
$\tilde g=H_{\tilde f} \circ \tilde f \circ H^{-1}_{\tilde f}$
belongs to $\cE^1_0$ for every $\tilde f \in \cU$.  
Consider the mapping $\Pi: \tilde f \mapsto \tilde g$; it is clearly continuous.
It is also open: for any $\hat{g} \in \cE^1_0$ close to $g$,
$H_f^{-1} \circ \hat{g}\circ H_f \in \cE^1$ is close to $f$ and 
is mapped by $\Pi$ to $\hat{g}$. 
The preimage of $\cR_0$ under $\Pi$ is thus a residual subset
of $\cU$, which contains only maps which have no aci$\sigma$.

\subsection{Lispchitz maps}

Let $\Elip$ be the set of Lipschitz
local homeomorphisms $f:\T^1 \to \T^1$, such that $f(0)=0$ and
$$
\lambda_f = \essinf_{x \in \T^1} |f'(x)|>1.
$$

We consider $\Elip$ endowed with the topology induced from
the Lipschitz metric:
$$
\dist(f,g)=\esssup_{x \in \T^1} |f'(x)-g'(x)|.
$$
We also let
$$
\Lambda_f=\esssup_{x \in \T^1} |f'(x)|.
$$

The distortion of the restriction of some iterate of $f$ to some interval is
$$
\Dist(f^n|J)=\frac {\esssup_{x \in J} |(f^n)'(x)|} {\essinf_{x \in J} |(f^n)'(x)|}\, .
$$
Clearly, if $f^n|J$ is 1-1 onto $I$ then 
for every measurable $X \subset J$ of positive measure,
$$
\frac{1}{\Dist(f^n|J)} \le
\frac{m(f^n(X))/m(I)}{m(X)/m(J)}
\le \Dist(f^n|J).
$$

%
%
%
%

\subsection{Recurrence properties}

We say that $f \in \Elip$ is \emph {ergodic} with respect to $m$ if every measurable
set $X$ such that $f^{-1}(X)=X$ satisfies $m(X)=0$ or $m(X)=1$.
We say that $f \in \Elip$ is \emph {conservative} with respect to $m$ if
every measurable set $X$ such that $X \cap \bigcup_{k=1}^\infty f^{-k}(X)=\emptyset$
satisfies $m(X)=0$.

It is easy to see that $f \in \Elip$ is ergodic and conservative with respect to $m$ if and
only if for every measurable set $X$ which is forward invariant (that is, $f(X) \subset X$)
we have $m(X)=0$ or $m(X)=1$.

If $f \in \Elip$ and $X \subset \T^1$ is a measurable set, 
then we denote $X_f$ the set of points in $X$ that return to $X$ by forward iteration by $f$.
It is easy to see that if $f$ is conservative then $m(X_f)=m(X)$.

We denote by $f_X : X_f \to X$ the first return map.

\begin{lemma}\label{l.preserve}
If $\mu$ is any $f$-invariant $\sigma$-finite measure then 
$\mu(f_X^{-1} Y) \le \mu(Y)$ for all measurable $Y \subset X$.
\end{lemma}
\begin{proof}
We can assume $\mu(Y) < \infty$.
For $n\ge 1$, let 
$$
Z_n = f^{-1} (X^c) \cap \cdots \cap f^{-(n-1)}(X^c) \cap f^{-n} (Y), \quad
Y_n = X \cap Z_n.
$$
Then $\bigsqcup_{n=1}^\infty Y_n = f_X^{-1}(Y)$.
Since $Z_{n+1} = f^{-1}(X^c \cap Z_n)$, we have
$$
\mu(Y_n) = \mu(X \cap Z_n) = \mu(Z_n) - \mu(X^c \cap Z_n) = \mu(Z_n) - \mu(Z_{n+1}).
$$
Therefore $\sum_{n=1}^\infty \mu(Y_n) \le \mu(Z_1) = \mu(f^{-1}(Y)) = \mu(Y)$.
\end{proof}

\subsection{Markov partitions}

Let $f \in \Elip$. 
The points in $f^{-n}(0)$ divide the circle into $|d_f|^n$ open intervals,
which are called \emph{Markov intervals of order $n$}.
The image of a Markov interval of order $n$ is a Markov interval of order $n-1$.
If $I$ is a Markov interval of order $n$ then 
$f^n| I$ is a 1-1 map onto $(0,1)$.

If $I$ is a Markov interval and $f_I$ is the first return map to $I$,
then there exist disjoint (Markov) intervals $I_j \subset I$ 
such that $f_I|I_j$ is an homeomorphism onto $I$ for each $j$.
If $f$ is conservative with respect to Lebesgue measure
then the intervals $I_j$ cover $I$ $m$-mod~$0$.

\subsection{Piecewise linear approximations}

Let us say that a map $f \in \Elip$ belongs to $\Epl$ 
if there exists $n \geq 1$ such that
for every Markov interval $I$ of order $n$, $f|I$ is linear.

\begin{lemma}\label{l.eplisdense}
The set $\Epl$ is dense in $\cE^1_0$ in the Lipschitz metric.
\end{lemma}

\begin{proof}
Given $f \in \cE_0^1$
and $n \ge 2$, we define a map $f_n: \T^1 \to \T^1$ as follows:
For each Markov interval $I$ of order $n$ for $f$, 
let $f_n$ map $I$ onto $f(I)$ linearly,
and so that $f_n$ equals $f$ in the boundary of $I$.
Clearly, $f_n \in \Epl$.
We claim that $f_n \to f$ in the Lipschitz metric.

Notice that the lengths of Markov intervals $I$ of order $n$ go uniformly to $0$ as $n \to \infty$:
in fact, $m(I) \leq \lambda_f^{-n}$.  
Since $f$ is $C^1$, for every $\delta>0$ we can choose $n_0$ such
that in each Markov interval $I$ of order $n \geq n_0$, $\sup_{x,y \in I} |f'(x)-f'(y)| \leq \delta$.
Then $\dist(f,f_n) \leq \delta$.
\end{proof}

\begin{lemma}\label{l.linearisergodic}
Every $f \in \Epl$ is  ergodic and conservative with respect to Lebesgue measure.
\end{lemma}

\begin{proof}
Given $f \in \Epl$, let $n$ be such $f$ is linear on Markov intervals of order $n$.
Let $X$ be a forward invariant set of positive Lebesgue measure.
Notice that if $I$ is a Markov interval of order $N \geq n$, then
$f^{N-n}|I$ is linear onto some Markov interval of order $n$, and $f^N|I$ is onto $(0,1)$
with distortion bounded by $\Lambda_f^n / \lambda_f^n$.  By the Lebesgue Density Points theorem,
for almost every $x \in X$, $\lim_{N \to \infty} m(I_N \cap X)/m(I_N)=1$, where $I_N$ is the Markov
interval of order $N$ containing $x$.  Applying $f^N$ and using the bound on the distortion, we see that
$m(f^N(X \cap I_N)) \to 1$ as $N \to \infty$.  So $m(X)=1$ and the result follows.
\end{proof}

\section{Plan of proof}

\begin{defn}\label{d.distorted}
Let $I$ be an interval and $\phi: I \to I$ non-singular (w.r.t.~the Lebesgue measure) map.
Then $\phi$ is called \emph{distorted} if 
there exists a measurable set $A \subset I$ such that
$\dfrac{m(A)}{m(I)}>.4$, and
$\dfrac{d(\phi_* m_I)}{dm_I}>2$ on $A$ (where $m_I$ is Lebesgue measure on $I$).
\end{defn}

Recall that
$$
\dfrac{d(\phi_* m_I)}{dm_I}(y) = \sum_{x \in f^{-1}(y)} \dfrac{1}{|\phi'(x)|}.
$$

\begin{defn}\label{d.good}
If $\delta>0$ and $f \in \Elip$ then
$f$ is called \emph{$\delta$-good} if there exists a 
family $\cI$ of Markov intervals (possibly of different orders), 
all of length at most $\delta$,
such that $m\left( \bigcup_{I \in \cI} I \right) > 1-\delta$
and for every $I \in \cI$, the first return map
$f_I: I \to I$ is distorted.
\end{defn}

\begin{figure}[htb]
\psfrag{A}{{\tiny $A$}} 
\begin{center}
\includegraphics[width=4cm]{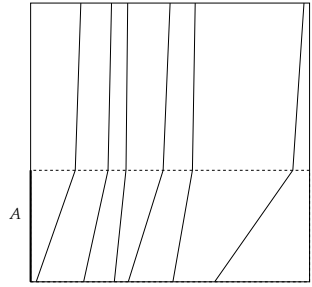}
\caption{A typical graph of a distorted $f_I$.
In this example, $A$ is an interval.
There should be infinitely many branches.}
\label{f.good}
\end{center}
\end{figure}

We can now state the three key technical results of this paper.

\begin{prop}\label{p.criterium}  
If $f \in \cE^1_0$ is $\delta$-good for every $\delta>0$ then $f$ has no aci$\sigma$.
\end{prop}


The next proposition says the condition of being $\delta$-good is open in two senses.

\begin{prop}\label{p.open}
Let $f\in\Elip$ be $\delta$-good for some $\delta>0$.
Then:
\begin{enumerate}
\item There exists $\beta >0$ such that if 
a map $\tilde f \in \Elip$ satisfies 
$$
m\left(\{x \in \T^1;\; \tilde f(x) \neq f(x)\} \right)< \beta
$$ 
then it is $\delta$-good.
\item Assuming that $f\in\cE^1_0$, there exists $\gamma>0$
such that if $\tilde f \in \cE^1_0$ satisfies $\dist(\tilde f, f) < \gamma$ then
$\tilde f$ is $\delta$-good.
\end{enumerate}
\end{prop}


\begin{prop}\label{p.dense}
For any $\delta>0$, the set of maps that are $\delta$-good
is dense in $\cE^1_0$ with the Lipschitz metric.
\end{prop}

Let us first see how to conclude theorem~\ref{t.reduced} (and hence, by the reduction,
theorem~\ref{t.main}) from the three propositions.

\begin{proof}[Proof of theorem~\ref{t.reduced}]
The set $\cU_\delta$ of $f \in \cE^1_0$ which are $\delta$-good 
is $C^1$-open, by part~(ii) of proposition~\ref{p.open},
and $C^1$-dense, by proposition~\ref{p.dense}. 
So $\cR_0=\bigcap_{\delta>0} \cU_\delta$ is a residual set of $\cE^1_0$, 
and by proposition~\ref{p.criterium} it consists of maps which do not have an aci$\sigma$.
\end{proof}

We now prove propositions~\ref{p.criterium} and \ref{p.open},
and leave the harder proof of proposition~\ref{p.dense}
(where part~(i) in proposition~\ref{p.open} is used)
for the next section.

\begin{proof}[Proof of proposition~\ref{p.criterium}]
Assume that for all $\delta>0$, $f \in \cE^1$ is $\delta$-good;
let $\cI_\delta$ be the corresponding family of intervals as in definition~\ref{d.good}.
Assume $f$ has an aci$\sigma$ $\mu$,
and let $\rho = \frac{d\mu}{dm}$ be its density.
Let $c>0$ be such that the set $Z = \{x \in \T^1; c \le \rho(x) \le 1.1 c \}$ has positive Lebesgue measure.
By a density point argument,
there exists an interval $I \in \bigcup_{\delta} \cI_\delta$ such that 
$m(I \cap Z) / m(I)> .9$.
Since $f_I$ is distorted, there exists $A \subset I$ such that 
$m(A)/m(I)=.4$ and $d(f_I)_* m /dm>2$ on $A$.
Let $Y = A \cap Z$.
We have
$$
m(Y) \ge m(A) - m(I \setminus Z) > .4 \, m(I) - .1 \, m(I) = .3 \,m(I)
$$
and
$$
\mu(Y) \le 1.1 c m(Y) \le 1.1 \times .4 c \, m(I) = .44 c \, m(I) 
$$
Moreover,
\begin{multline*}
\mu \left(f_I^{-1}(Y)\right)        \ge 
c \, m\left(Z \cap f_I^{-1}(Y)\right)  \ge
c \left[ m\left(f_I^{-1}(Y)\right) - m(I \setminus Z) \right] \\ >
c [2 \, m(Y) - .1 \, m(I)] \ge
c [2 \times .3 - .1] m(I) =
.5 c \, m(I)
> \mu(Y).
\end{multline*}
This contradicts lemma~\ref{l.preserve}.
\end{proof}

\begin{proof}[Proof of proposition~\ref{p.open}]
Let $\cI$ be the family of Markov intervals as in definition~\ref{d.good};
clearly we can assume it is finite,
say, $\cI = \{I_i; \; 1 \le i \le i_0\}$.
Let $A_i \subset I_i$ be the set that 
gets enlarged under $(f_{I_i})^{-1}$ according to definition~\ref{d.distorted}.
Let $J_{i,1}$, $J_{i,2}$, \ldots be the connected components of the domain of $f_{I_i}$,
and let $n_{i,j}$ be such that $f^{n_{i,j}} (J_{i,j}) = I_i$.
Then 
$$
\sum_{j=1}^\infty |((f^{n_{i,j}}|J_{i,j})^{-1})'(y)| > 2
\quad\text{for every $y\in A_i$.}
$$
Slightly reducing the sets $A_i$ (still keeping $m(A_i)/m(I_i)>.4$), 
we can find $j_0$ such that
$$
\sum_{i=1}^{j_0} |((f^{n_{i,j}}|J_{i,j})^{-1})'(y)| > 2 + \epsilon
\quad\text{for every $y\in A_i$, for every $i=1,\ldots,i_0$,}
$$
where $\epsilon$ is some fixed positive number.
Also let $N = \max \{ n_{i,j}; \; 1\le i \le i_0, \ 1 \le j \le j_0 \}$.

If $\tilde f$ is another map in $\Elip$ which is $C^0$-close to $f$,
then for each interval $I_i \in \cI$ there is an interval 
$\tilde I_i$ which is Markov for $\tilde f$ and is close to $I_i$.
Clearly if the $C^0$-distance between $\tilde f$ and $f$ is sufficiently small
then each $\tilde I_i$ has length $<\delta$, and their union has measure $>1-\delta$.
Further, for each $\tilde I_i$ there exist intervals
$\tilde J_{i,j}$, $1 \le j \le j_0$, which are close to $J_{i,j}$
and such that $\tilde f^{n_{i,j}}(\tilde J_{i,j}) = \tilde I_i$.

With these notations fixed, we complete
the proofs of the two parts of the proposition  separately.

\smallskip

\emph{Part~(i):} 
Let $\tilde f \in \Elip$ so that the set 
$U = \{x \in \T^1;\;\tilde f(x) \neq f(x)\}$ has 
$m(U) < \beta$,
where how small $\beta$ needs to be will become clear along the way.
First, notice that the $C^0$-distance between $\tilde f$ and $f$ is small
(in fact, less than $\Lambda_f \beta$).
So we can define intervals 
$\tilde I_i$, $\tilde J_{i,j}$, for $1 \le i \le i_0$ and $1 \le j \le j_0$
as explained above.
Let
$V = \bigcup_{n=0}^N f^n(U)$.
Then $m(V) \le (1+\Lambda_f+ \cdots+\Lambda_f^N)\beta$ is small.
Define
$\tilde A_i = \tilde I_i \cap A_i \setminus V$.
Then $m(A_i \setminus \tilde A_i)$ is small: 
at most $m(I_i \setminus \tilde I_i) + m(V)$.
So if $\beta$ is small enough then $m(\tilde A_i) / m(\tilde I_i) > .4$.
Moreover, if $y \in \tilde A_i$ then
$\dfrac{d(\tilde f_{\tilde I_i})_* m}{dm} (y) > 2 + \epsilon$.
This shows that $\tilde f_{\tilde I_i}$ is distorted for each $i$
and accordingly that $\tilde f$ is $\delta$-good.

\smallskip

\emph{Part~(ii):} 
Now assume $f$ is $C^1$ and $\tilde f$ is $\gamma$-$C^1$-close to $f$.
Again we can define intervals 
$\tilde I_i$, $\tilde J_{i,j}$, for $1 \le i \le i_0$ and $1 \le j \le j_0$.
Let $\tilde A_i = \tilde I_i \cap A_i$. 
By taking a small $\gamma$,
we guarantee that $\tilde f^n$ is $C^1$-close to $f^n$ for $1\le n \le N$,
and therefore
$$
\sum_{i=1}^{j_0} |((\tilde f^{n_{i,j}}| \tilde J_{i,j})^{-1})'(y)| > 2
\quad\text{for every $y\in \tilde A_i$.}
$$
So the $\tilde f_{\tilde I_i}$ are distorted and
and $\tilde f$ is $\delta$-good.
\end{proof}

We remark that with a little more effort it is possible to improve 
simultaneously the two parts of proposition~\ref{p.open}, 
showing that being $\delta$-good is an open condition in $\Elip$ 
in the \emph{bounded variation metric}
$d_\mathrm{BV}(f,g) = \int_{\T^1}|f'-g'| \, dm$.

\section{Proof of proposition~\ref{p.dense}}

Let $f_0 \in \cE_0^1$ and $\delta>0$;
we will show that there exists a $\delta$-good map $h \in \cE_0^1$ such that
$\dist(h, f_0) < 3\delta$.
For simplicity, we will assume that $f$ is orientation-preserving.  
The proof can be easily adapted to cover the general case.

\subsection{Step 1. Linearization}
By lemma~\ref{l.eplisdense}, we can find $f \in \Epl$
with $\dist(f, f_0)<\delta$.
Since $f \in \Epl$, there exists $n_0$ such that 
if $I$ is a Markov interval of order $n \ge n_0$ then $f|I$ is linear.

Let $\ell \geq n_0$ be large (to be specified later).
Fix a Markov interval $T$ such that the sets
$T$, $f^{-1}(T)$, \dots, $f^{-\ell}(T)$ are disjoint,
and their union has Lebesgue measure less than $\delta$.

Let $\cP_T$ be the collection of (Markov) subintervals of $T$ that are sent 
onto $T$ by $f_T$.  
Let $K_T=T \setminus \bigcup_{I \in \cP_T} I$.
Notice that for any $I \in \cP_T$, $\order(I) \geq \order(T)+\ell$,
$\order(T) \geq \ell \geq n_0$ and
$f_T | I = f^{\order(I) - \order(T)} |I$ is linear.

\subsection{Step 2. Another perturbation}

If $I$ is an interval, 
denote by $\Phi_I$ the only order-preserving linear bijection $I\to(0,1)$.

Each interval in $\cP_T$ has length at most $\lambda_f^{-\ell} m(T)$.
Since $\ell$ is large, we can find 
$\xi \in T \setminus \bigcup_{I \in \cP_T} I$ such that
$$
\eta=\Phi_T(\xi) \in (.4, .41). 
$$

\begin{figure}[htb]
\psfrag{etaemeio}{{\tiny $\eta + \tfrac{1}{2}$}} 
\psfrag{eta}{{\tiny $\eta$}} 
\begin{center}
\includegraphics[width=4cm]{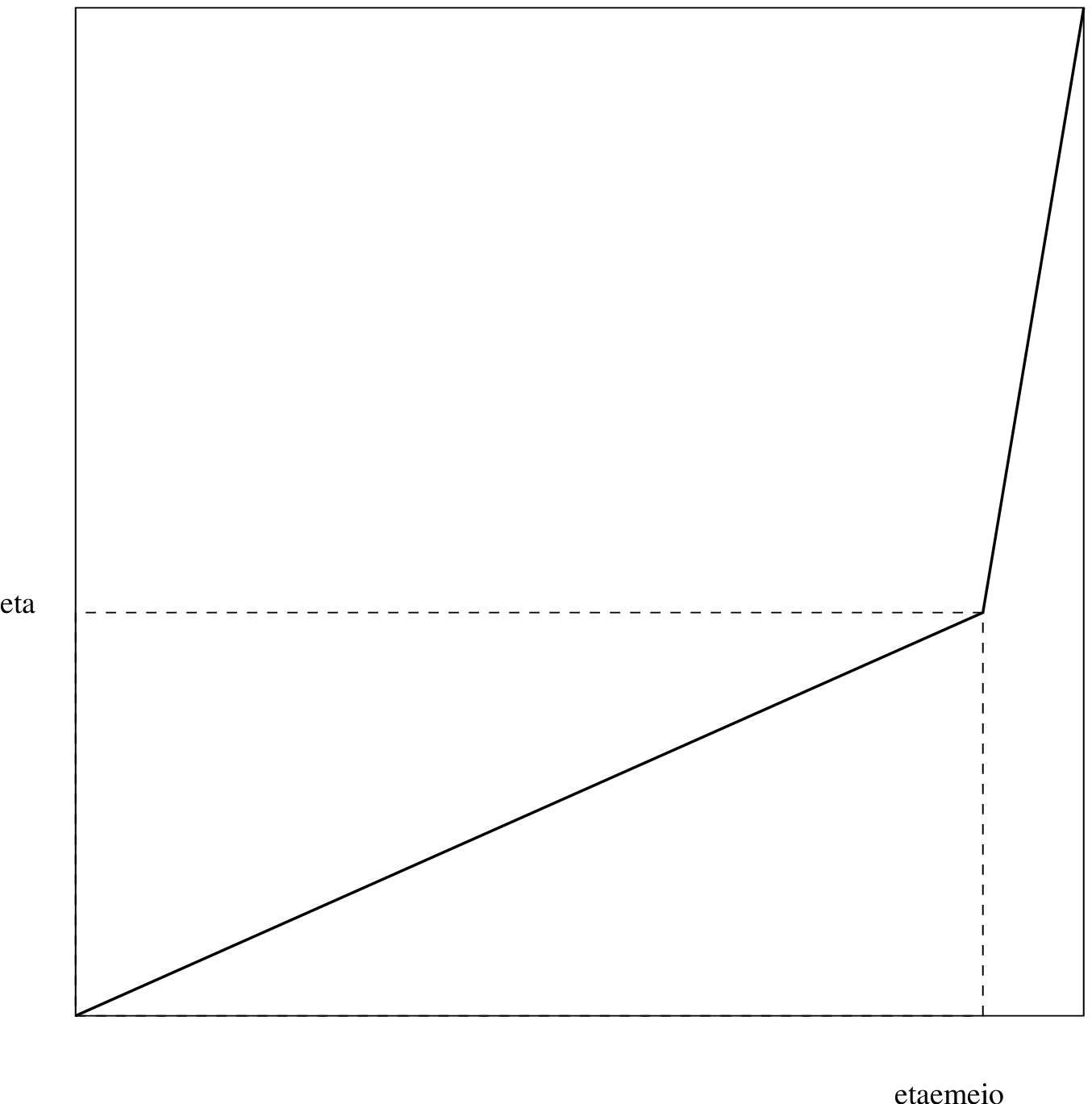}
\caption{Graph of $Q:(0,1) \to (0,1)$.}
\label{f.q}
\end{center}
\end{figure}

We will define a perturbation $g$ of $f$ as follows:
\begin{itemize}
\item $g$ equals $f$ outside $f^{-1}(T \setminus K_T) \cup \cdots \cup f^{-\ell}(T \setminus K_T)$.

\item Let $\Xi=f_T^{-1}(\xi)$.
Consider all sequences of intervals
$$
I_\ell \to I_{\ell-1} \to \cdots \to I_0
$$
with $f(I_j) = I_{j-1}$ and $I_0 \in \cP_T$.
Let $\xi_0=I_0 \cap \Xi$.
Then $\Phi_{I_0}(\xi_0)=\eta$.
Let $\xi_i=\Phi^{-1}_{I_i}(\eta+\frac {i} {2 \ell})$.  We define $g|I_i$ for $1 \leq i \leq \ell$ as the unique orientation preserving homeomorphism onto $I_{i-1}$ whose restriction to each connected component of $I_i \setminus \{\xi_i\}$ is linear and such that $g(\xi_i)=\xi_{i-1}$.  Let $Q=\Phi_{I_0} \circ g^\ell \circ \Phi^{-1}_{I_\ell}$.
It is the homeomorphism of $(0,1)$ depicted in figure~\ref{f.q}.
\end{itemize}

Notice that the Lipschitz distance between $g$ and $f$ is at most $C_f/\ell$, 
for some constant $C_f$ depending on $f$ only, and hence it is $<\delta$ 
since $\ell$ is large.

\subsection{Step 3. Some properties of $g$}

Let $\cP^k_{g,T}$ be the 
family of ($g$-Markov) subintervals of $T$ that are sent onto $T$ by $g_T^k$. 
Notice that $\cP^1_{g,T}=\cP_T$.

If $\zeta$ is a singularity for $g$, 
in the sense that $g$ is not linear in a neighborhood of $\zeta$,
then:
\begin{enumerate}
\item either $\zeta$ is a singularity for $f$ -- in this case it is contained in $f^{-n_0}(0)$;
\item or $\zeta$ belongs to $f^{-i}(K_T)$ for some $1 \leq i \leq \ell$;
\item or $\zeta$ belongs to $g^{-i}(\Xi)$ for some $1 \leq i \leq \ell$.
\end{enumerate}

\begin{lemma} \label{l.g}
If $L$ is an element of $\cP^2_{g,T}$ and
$J$ is a connected component of $g^{-k}(L)$, $k \ge 1$,
then $g^k | J$ is linear.
\end{lemma}

\begin{proof}
If $g^k|J$ is not linear then $g^j(J)$ intersects 
a singularity $\zeta$ of $g$ for some $0 \leq j<k$.  
Since $L = g^k(J)$ belongs to $\cP^2_{g,T}$, there are at least three different $i>0$
such that $g^i(\zeta) \in T$.  
Thus $\zeta$ cannot be of the type~(i) above:
indeed, the set $f^{-n_0}(0)= g^{-n_0}(0)$ is forward invariant 
for both $f$ and $g$, and it does not intersect $T$.
$\zeta$ cannot be of the type~(ii): 
the first iterate of $\zeta$ that belongs to $T$ belongs indeed
to $K_T$, and subsequent iterates do not enter $T$ again.  
$\zeta$ cannot be a singularity of type~(iii):
the first iterate of $\zeta$ that belongs to $T$ also belongs to $\Xi$, 
so the second iterate that belongs to $T$ is $\xi$, 
and the subsequent iterates lie outside $T$.  
So there can be no such singularity, and the result follows.
\end{proof}

\begin{lemma} \label{l.distortion}
If $L$ is an element of $\cP^k_{g,T}$ then $\Dist(g_T^k|L) \leq
\left ( \Lambda_g / \lambda_g \right )^{2 \ell}$.
\end{lemma}

\begin{proof}
Take $L \in \cP^1_{g,T}$.
Let $r$ be such that $g^r(L)=T$.  
Then $g^{r-\ell}|L$ is linear, and hence 
$\Dist(g_T|L)=\Dist(g^\ell|g^{r-\ell}(L)) \leq \left (\Lambda_g / \lambda_g \right )^\ell$.
This implies the assertion of the lemma for $k=1$ and $k=2$.
Now, if $k \ge 2$ and and $L \in \cP^k_{g,T}$ then,
by lemma~\ref{l.g}, $g_T^{k-2}|L$ is linear,
so the assertion also follows.
\end{proof}

\begin{lemma}\label{l.ergodic}
$g$ is ergodic and conservative with respect to Lebesgue measure.
\end{lemma}

\begin{proof} 
We will adapt the argument in the proof of lemma~\ref{l.linearisergodic}.
Let $X \subset \T^1$ be a forward $g$-invariant set with $m(X)>0$.

Assume that $X \cap T$ has zero Lebesgue measure. 
By a density point argument, we can take a $g$-Markov interval $L$ 
such that $m(L \cap X)/m(L)$ is close to $1$.
Any forward-image $g^k(L)$ cannot be contained 
in the set $W = \bigcup_{j=1}^\ell f^{-j}(T)$.
By the Markov property, $g^k(L) \cap W = \emptyset$ for $0 \le k \le \order(L)-\order(T)-\ell$.
Since $g$ equals $f$ outside $W$,
$g^{\order(L)-\order(T)-\ell}|L$ is linear.
In particular, 
$$
\Dist(g^{\order(L)}|L) \leq \left ( \Lambda_g / \lambda_g \right )^{\order(T)+\ell}.
$$
It follows that $m(g^{\order(L)}(L \cap X))$ is close to $1$.
This shows that the assumption $m(X \cap T) = 0$ cannot be true.

Since $f$ is conservative (lemma~\ref{l.linearisergodic}), 
$m$-almost every point in $T$ returns to $T$ 
by forward iterates of $f$.
It follows that the same is true for $g$.
So the intervals in $\cP^k_{g,T}$ cover $T$ $m$-mod $0$. 
Moreover, these intervals have lengths at most $\lambda_g^{k \ell}$.
Therefore we can find $L$ and $k \geq 2$ such that $L \in \cP^k_{g,T}$ and 
$m(L \cap X)/m(L)$ is arbitrarily close to $1$.  By lemma~\ref{l.distortion},
$m(g_T^k(L \cap X))/m(T)$ is arbitrarily close to $1$.  
It follows that $m(X \cap T)/m(T)=1$.  
So $g^{\order(T)}(X \cap T) \subset X$ has full $m$-measure on the circle.
\end{proof}

\begin{lemma}
$g$ is $\delta$-good.
\end{lemma}

\begin{proof}
First, let us define the family $\cI$:
an interval $I$ belongs to $\cI$ iff
there exists $n=n(I) > \ell$ such that $f^n(I) = T$ and $f^k(T) \cap T=\emptyset$ for $0 \leq k<n$.
Notice that:
\begin{itemize}
\item For every $I \in \cI$, $g^k(I)=f^k(I)$ for $0 \leq k \leq n(I)$, and hence $I$ is Markov for $g$.
\item For every $I \in \cI$, $m(I)<m(T)<\delta$.
\item
$\bigcup_{I \in \cI} I  = \T^1 \setminus \bigcup_{j=0}^\ell f^{-j}(T)$ $m$-mod $0$;
in particular $m(\bigcup_{I \in \cI} I) > 1-\delta$.
\end{itemize}

We have to show that for each $I \in \cI$, $g_I$ is distorted.  
For this, it is enough to prove that 
(where $\cP_{g,I}$ is the collection of intervals $J \subset I$ that are sent onto $I$ by $g_I$):
\begin{enumerate}
\item $\bigcup_{J \in \cP_{g,I}} J=I$ $m$-mod $0$;
\item if $J \in \cP_{g,I}$ and $r$ is such that $g^r(J)=I$ then $\Phi_I \circ g^r \circ \Phi^{-1}_J=Q$.
\end{enumerate}
Indeed in this case we can take $A=\Phi_I^{-1}(0,\eta)$ in definition~\ref{d.distorted}.

The first property follows from lemma~\ref{l.ergodic}.
Let us check the second one.
Let $\{i_1<\cdots<i_t\}=\{k; \; 0 \leq k \leq r, \  g^k(J) \subset g^{-\ell}(T)\}$.  
Then 
$$
0<i_1<i_1+\ell<\cdots <i_{t-1}<i_{t-1}+\ell<i_t<i_t+\ell<r.
$$

We claim that $g^{i_t}|J$ is linear.  This is clear if $t=1$.
Notice that $g^{i_t+\ell}(J)$ is an element of $\cP_T$: 
indeed, $g^{r-i_t-\ell+n(I)}$ takes $g^{i_t+\ell}(J)$
onto $T$.  This implies that $g^{i_j+\ell}(J)$ is an element of $\cP^{1+t-j}_{g,T}$.  
By lemma~\ref{l.g},
$g^{i_{t-1}+\ell}|J$ is linear, therefore $g^{i_t}|J$ is linear, as claimed.

Since $g^{i_t+\ell}(J)$ is an element of $\cP_T$, $\Phi_{g^{i_t+\ell}(J)} \circ g^\ell \circ \Phi^{-1}_{g^{i_t}(J)}=Q$.
It is also clear that $g^{r-i_t-\ell}|g^{i_t+\ell}(J)$ is linear.  It follows that $\Phi_I \circ g^r \circ \Phi^{-1}_J=Q$, as desired.
\end{proof}

\subsection{Step 4. Smoothening $g$} 

For $X \subset \T^1$ and $\epsilon>0$,
let $B_\epsilon(X)$ be the $\epsilon$-neighborhood of $X$.
Let $S$ be the (already described) set of singularities of $g$;
then $S$ is closed in $\T^1$ and $m(S)=0$.  
In particular, $m(B_\alpha(S)) \to 0$ as $\alpha \to 0$.

Let $G:\R \to \R$ be the lift of $g$ satisfying $G(0)=0$.
For $\alpha>0$, let 
$$
G_\alpha(x) = \frac {1} {2 \alpha} \int_{x-\alpha}^{x+\alpha} G(y) dy.
$$  
Then $G_\alpha:\R \to \R$ is the lift of some $C^1$ map
$g_\alpha:\T^1 \to \T^1$ such that for every $x \in \T^1$,
$$
|g_\alpha'(x)-f'_0(x)| \leq \dist(g, f_0)+\sup_{y \in B_\alpha(x)} |f'_0(y)-f'_0(x)|.
$$
(Recall $f_0$ is the original $C^1$ map.)
Hence $|g'_\alpha(x)-f'_0(x)| \le 2\delta$ for every $x \in \T^1$ if $\alpha$ is sufficiently small.
Moreover, since $g$ is linear on connected components of $\T^1\setminus S$,
$g_\alpha$ equals $g$ outside the $\alpha$-neighborhood of $S$.

Unfortunately, $g_\alpha$ does not necessarily fix $0$.
To remedy that, take a family of $C^1$ diffeomorphisms 
$\phi_\xi:\T^1 \to \T^1$ parameterized by $\xi \in \T^1$, 
such that 
$$
\phi_\xi(\xi)=0, \quad 
\lim_{\xi \to 0} \sup_{x \in \T^1} |\phi'_\xi(x)-1| = 0, \quad\text{and}\quad 
\lim_{\xi \to 0} m\{x \in \T^1;\; \phi_\xi(x) \neq x\} = 0.
$$
Define $h_\alpha=\phi_{g_\alpha(0)} \circ g_\alpha$.  
If $\alpha$ is small then
$h_\alpha \in \cE^1_0$ and $\dist(h_\alpha, f_0) < 3\delta$.
Also, $\lim_{\alpha \to 0} m\{x \in \T^1;\; h_\alpha(x) \neq g(x)\}=0$.
By part~(i) of proposition~\ref{p.open}, $h_\alpha$ is $\delta$-good 
provided $\alpha$ is small enough.
This concludes the proof of proposition~\ref{p.dense}.


\bigskip

\begin{ack}
This research was partially conducted during the period A.A. served as a Clay Research Fellow.
The paper was written while J.B.~was visiting IMPA, partially supported by Procad/Capes project.
J.B.~thanks Philippe Thieullen for providing several references and for
stimulating conversations.
\end{ack}


\end{document}